\documentclass{amsart}

\newtheorem{theorem}{Theorem}[section]
\newtheorem{lemma}[theorem]{Lemma}

\theoremstyle{definition}
\newtheorem{definition}[theorem]{Definition}
\newtheorem{example}[theorem]{Example}
\newtheorem{xca}[theorem]{Exercise}

\theoremstyle{remark}
\newtheorem{remark}[theorem]{Remark}

\numberwithin{equation}{section}

\newcommand{\abs}[1]{\lvert#1\rvert}

\newcommand{\blankbox}[2]{%
  \parbox{\columnwidth}{\centering
    \setlength{\fboxsep}{0pt}%
    \fbox{\raisebox{0pt}[#2]{\hspace{#1}}}%
  }%
}

\author{Gizem Karaali}
\address{Department of Mathematics, University of California, Santa Barbara, Ca 93106}
\email{gizem@math.ucsb.edu}



\newtheorem{prop}{Proposition}

\theoremstyle{definition}
\newtheorem*{acknow}{Acknowledgments}

\newcommand{\baseRing}[1]{\ensuremath{\mathbb{#1}}}
\newcommand{\C}{\baseRing{C}}

\newcommand{\Z}{\baseRing{Z}}

\newcommand{\g}{\mathfrak{g}}
\newcommand{\p}{\mathfrak{p}}
\renewcommand{\d}{\mathfrak{d}}
\renewcommand{\t}{\mathfrak{t}}
\newcommand{\s}{\mathfrak{s}}

\renewcommand{\phi}{\varphi}

\renewcommand{\Im}{\operatorname{Im}}

\begin{document}

\title{A New Lie Bialgebra Structure on $sl(2,1)$}

\begin{abstract}
We describe the Lie bialgebra structure on the Lie superalgebra $sl(2,1)$ related to an $r-$matrix
that cannot be obtained by a Belavin-Drinfeld type construction. This structure makes $sl(2,1)$
into the Drinfeld double of a four-dimensional subalgebra.
\end{abstract}

\maketitle

It is well-known that non-degenerate $r-$matrices (describing quasitriangular Lie bialgebra
structures) on simple Lie algebras are classified by Belavin-Drinfeld triples, (the original
references are \cite{BD1, BD2}, more pedagogical presentations providing ample background can be
found in \cite{CP, ES}). A similar construction using Belavin-Drinfeld type triples is possible for
simple Lie superalgebras with nondegenerate Killing form. Surprisingly, though, in the super
setting, there are certain non-degenerate $r-$matrices that do not fit such a description, see
\cite{Kar}. 

The purpose of this note is to study the super Lie bialgebra structure associated to such an
$r-$matrix on the simple Lie superalgebra $sl(2,1)$. We start in Section
\ref{SectionSuperLieBialgebras} with some background on super Lie bialgebra structures and some
basic constructions related to them. In Section \ref{SectionExamplermatrix}, we explicitly
describe the $r$-matrix that we will be interested in. Section \ref{SectionLieBialgebraStructure}
describes in detail the associated super Lie bialgebra structure on $sl(2,1)$; this structure makes
$sl(2,1)$ into the Drinfeld double of a four-dimensional subalgebra. A comparison with the standard
super Lie bialgebra structure is also provided in this section. We end in Section
\ref{SectionConclusion} with a brief discussion of the results and further directions for
investigation.

\begin{acknow}
The author would like to thank N. Reshetikhin, V. Serganova and M. Yakimov for their comments and
suggestions. 
\end{acknow}

\section{Lie Bialgebra Structures on Lie Superalgebras}
\label{SectionSuperLieBialgebras}

\subsection{Cohomology of Lie superalgebras}
\label{SuperCohomology}

The cohomology theory of Lie superalgebras is more complicated than that of Lie algebras. Even
for simple Lie superalgebras and for low dimensions, it is not yet completed. Here we summarize
certain basic facts that we will use. For more on the cohomology theory of Lie superalgebras
one can look at \cite{Fu,ScZh}.

Recall that if $\g$ is a Lie algebra, then an $n-$cochain taking values in a $\g-$module $M$ is an
alternating $n-$linear map $f(x_1, x_2, \cdots, x_n)$ of $n$ variables in $\g$. We can view each
such $n-$cochain as a linear map ${f : \bigwedge^n\g \rightarrow M}$. In this case, the coboundary
$df$ of an $n-$cochain $f$ is the ${(n+1)-}$cochain defined by: 
\begin{eqnarray*} 
df(x_1, \cdots, x_{n+1}) 
&=& \sum_{i=1}^{n+1} (-1)^{i+1} x_i f(x_1, \cdots, \hat{x_i},
\cdots, x_{n+1}) \\  
&+& \sum_{1 \le i<j \le n+1} (-1)^{i+j} f( [x_i, x_j] , x_1,
\cdots , \hat{x_i} , \cdots , \hat{x_j} , \cdots , x_{n+1})
\end{eqnarray*}

If $\g$ is a Lie superalgebra, then the space of $n-$cochains with values in a $\g-$module $M =
{M_{\overline{0}} \oplus M_{\overline{1}}}$ is itself a graded space. Denoting this space by
$C^n(\g, M)$, we have: 
\[ C^n(\g, M) = \bigoplus_{i+j = n} Hom(\wedge^i
\g_{\overline{0}} \otimes S^j \g_{\overline{1}}, M).\]
\noindent 
The even part of $C^n(\g,M)$ is: 
\[ C^n_{\overline{0}}(\g,M) = \bigoplus_{i+j = n} Hom(\wedge^i
\g_{\overline{0}} \otimes S^j \g_{\overline{1}},
M_{\overline{j}}),\] 
\noindent 
while the odd part is given by:
\[ C^n_{\overline{1}}(\g,M) = \bigoplus_{i+j = n} Hom(\wedge^i
\g_{\overline{0}} \otimes S^j \g_{\overline{1}},
M_{\overline{j+1}}).\]

Equivalently we can view $C^n(\g,M)$, for $n \ge 1$, as the $\Z / {2\Z}$-graded vector space of all
super alternating $n-$linear maps $f$ of $\g^n = \g \times \g \times \cdots \times \g$ into $M$,
i.e. maps $f$ satisfying:  
\[ f(x_1, x_2, \cdots, x_i, x_{i+1}, \cdots, x_n) =
-(-1)^{|x_i||x_{i+1}|}
f(x_1, x_2, \cdots, x_{i+1}, x_i, \cdots, x_n). \]

We set $C^0(\g,M) = M$. The differential $d$ is defined as follows: For an $n-$cochain $f$, the
coboundary $df$ is an $(n+1)-$cochain given by:
\begin{align*}
df(x_1, \cdots, x_{n+1}) \quad 
&= \quad \sum_{i=1}^{n+1} \sigma_1(i,j) x_i f(x_1, \cdots, \hat{x_i},
\cdots, x_{n+1}) \\  
&+\quad  \sum_{1 \le i<j \le n+1} \sigma_2(i,j) f( [x_i,
x_j] , x_1, \cdots , 
\hat{x_i} , \cdots , \hat{x_j} , \cdots ,
x_{n+1})   
\end{align*}
\noindent
where the signs in the above sums are as follows:
\begin{align*}
\sigma_1(i,j) &= (-1)^{i+1}(-1)^{|x_i|(|f|+|x_1|+|x_2| +
\cdots + |x_{i-1}|)}\\
\sigma_2(i,j) &= (-1)^{i+j}(-1)^{|x_i||x_j|} 
(-1)^{|x_i|(|x_1| + \cdots + |x_{i-1}|)}
(-1)^{|x_j|(|x_1| + \cdots + |x_{j-1}|)} 
\end{align*}
\noindent
We note that this formula for $d$ agrees with that of \cite{ScZh} when we use the super alternating
property of $f$.

If $M = \g \otimes \g$, then $\g$ acts on $M$ on the left by the following extension of the adjoint
representation:
\[ g \cdot (a \otimes b) = (g \cdot a) \otimes b + (-1)^{|g||a|}a
\otimes (g \cdot b) = [g, a] \otimes b + (-1)^{|g||a|}a \otimes
[g, b] \] 
\noindent
In this setup, a $0-$cochain is a linear map ${f_0 : \C \rightarrow \g \otimes \g}$. Therefore it is
determined uniquely by ${f_0(1) \in \g \otimes \g}$ and hence can be identified with an element
$r$ of ${\g \otimes \g}$. The coboundary $dr$ of $r$ is a $1-$cochain defined by:
\[ dr(a) = a \cdot r = [a \otimes 1 + 1 \otimes a, r].\]
\noindent
A $1-$cochain is a linear map ${f : \g \rightarrow \g \otimes \g}$. It is a $1-$cocycle if ${df =
0}$, or in other words:
\begin{eqnarray*}
0 &=& df(a,b) = (-1)^{|a||f|}a \cdot f(b) - 
(-1)^{|b|(|f|+|a|)}b \cdot f(a) - f([a,b]) \\ 
&=& (-1)^{|a||f|} [a \otimes 1 + 1 \otimes a , f(b)] -
(-1)^{|b|(|f|+|a|)}
[b\otimes 1 + 1 \otimes b, f(a)] - f([a,b]) 
\end{eqnarray*}
\noindent
which we can rewrite as:
\begin{eqnarray*}
f([a,b]) &=& (-1)^{|a||f|} [a \otimes 1 + 1 \otimes a , f(b)] -
(-1)^{|b|(|f|+|a|)}
[b\otimes 1 + 1 \otimes b, f(a)] \\
&=&
[f(a), b \otimes 1 + 1 \otimes b] + (-1)^{|a||f|}[a \otimes
1 + 1 \otimes a, f(b)] \\
&=& (-1)^{|a||f|}a \cdot f(b) - (-1)^{|b||f(a)|}b \cdot f(a).
\end{eqnarray*}
\noindent
We will call the resulting formula the \emph{super cocycle condition}:
\[ f([a,b]) = (-1)^{|a||f|}a \cdot f(b) - (-1)^{|b||f(a)|}b
\cdot f(a).\]

\subsection{Super Lie bialgebras}
\label{SuperLieBialgebras}

A \emph{super Lie bialgebra} is a triple $(\g, [\cdot \; , \cdot ], \delta)$ such that:

\begin{enumerate}

\item $\g$ is a Lie superalgebra with the super bracket $[\cdot \;,\cdot]$;

\item $\delta : \g \rightarrow \g \otimes \g$ is a
skew-symmetric linear map whose dual ${\delta^* : \g^* \otimes
\g^* \rightarrow \g^*}$ defines a Lie superalgebra structure on
$\g^*$;

\item $\delta$ and $[\cdot \; , \cdot]$ are compatible in the following
sense:
\[ \delta([a,b]) = [\delta(a), b \otimes 1 + 1 \otimes b] + [a
\otimes 1 + 1 \otimes a, \delta(b)].\]
\end{enumerate}
\noindent
We will denote such a super Lie bialgebra by ${(\g,\delta)}$ if the super bracket
${[\cdot\;,\cdot]}$ is unambiguous. Note that the last condition is equivalent to $\delta$ being a
$1-$cocycle on $\g$ with values in $\g \otimes \g$, for the cohomology theory of Lie superalgebras
as summarized in \ref{SuperCohomology}. Since $\delta^*$ is a super bracket, $\delta$ is even, and
the super cocycle condition above coincides with the non-graded version.

The Jacobi identity for $\delta^*$ is equivalent to the following \emph{coJacobi identity} for
$\delta$ which holds for any $x \in \g$:
\[ 
Alt_s(\delta \otimes Id) \cdot \delta(x) = 0.\]
\noindent
Here ${Alt_s : \g \otimes \g \otimes \g \rightarrow \g \otimes
\g \otimes \g}$ is defined on homogeneous basis vectors by:
\[ Alt_s(a \otimes b \otimes c) =
a \otimes b \otimes c + (-1)^{|a|(|b|+|c|)}b \otimes c \otimes
a + (-1)^{|c|(|a|+|b|)}c \otimes a \otimes b.\] 

A (\emph{finite dimensional}) \emph{super Manin triple} is a triple  ${(\g, \g_+, \g_-)}$ of
(finite dimensional) Lie superalgebras such that: 
\begin{enumerate}
\item $\g$ is equipped with a non-degenerate super-symmetric invariant bilinear form
${(\cdot,\cdot)}$;
\item $\g_+$ and $\g_-$ are Lie subsuperalgebras of $\g$ and $\g = {\g_+ \oplus \g_-}$ as vector
spaces;
\item $\g_+$ and $\g_-$ are isotropic with respect to ${(\cdot,\cdot)}.$
\end{enumerate}
\noindent
Since the bilinear form is non-degenerate, $\g_+$ and $\g_-$ are in fact maximal isotropic or
\emph{Lagrangian} subsuperalgebras. We note here that in the infinite dimensional case, the
definition of a super Manin triple will be the same as above. However, one needs to take into
account the topology on vector spaces.

These two notions (i.e. super Lie bialgebras and super Manin triples) are related to one another in
a way similar to the Lie algebra case:

\begin{prop}
Let ${(\p, [\;,\;],\delta)}$ be a super Lie bialgebra and set $\g_+ = \p$ and $\g_- = \p^*$ and
define $\g = \g_+ \oplus \g_-$. Then ${(\g, \g_+, \g_-)}$ is a super Manin triple. Conversely,
any finite dimensional super Manin triple ${(\g,\g_+,\g_-)}$ gives rise to a super Lie bialgebra
structure on $\g_+$.
\end{prop}

\begin{remark}
This is Proposition $1$ of \cite{And} where it was proven modulo certain calculations left to the
reader. 
\end{remark}

\subsection{The Drinfeld double construction}
\label{DrinfeldDouble} 

Another related construction is that of the Drinfeld double. Here we will use a direct analogue
of the non-graded version, in the spirit of \cite{Gou}. Before explicitly presenting this approach,
we should also mention that other superizations of the double
construction exist. See for instance \cite{Vor} for a more geometrically motivated development of
the double.

Let $(\g, [\cdot \; , \cdot ]_{\g}, \delta_{\g})$ be a finite dimensional super Lie
bialgebra\footnotemark. Then clearly $\g^*$ is also a super Lie bialgebra with the associated
structures defined by: 
\[ [\cdot , \cdot ]_{\g^*} = (\delta_{\g})^* \qquad \delta_{\g^*} = ([\cdot \; , \cdot ]_{\g})^*.\]

\footnotetext{We will be assuming finite dimensionality, as this will be sufficient for our
purposes. Infinite dimensional analogues will be more technically involved, and since we do not
need them here, will not be discussed any further.}

Let us fix a homogeneous basis $\{ e_i \}$ for $\g$ and define the structure constants $C_{k}^{ij}$,
$D_{k}^{ij}$ of the relevant structures on $\g$ as follows:
\[ [e_i,e_j]_{\g} = \sum_{k} C^{k}_{ij}e_k \qquad \delta_{\g} (e_k) = \sum_{i,j} D_{k}^{ij} e_i 
\wedge e_j\]
\noindent
where we use the notation: $a \wedge b = a \otimes b - (-1)^{|a||b|}b \otimes a$ for any two
homogeneous elements $a$, $b$. 
From these we can determine the structure constants of $\g^*$; if we let $\{ e_i^*\}$ be the
homogeneous basis for $\g^*$ dual to $\{e_i\}$, then we have:
\[ [e_i^*,e_j^*]_{\g^*} = \sum_{k} C^{ij}_{k} e_k^* \qquad \delta_{\g^*} (e_k^*) = 
\sum_{i,j}  D^k_{ij} e_i^*  \wedge e_j^*\]
\noindent
where:
\[ C^{ij}_k = \left \{ \begin{matrix}
(-1)^{|e_i||e_j|} D_{k}^{ij} & i \neq j \\ \\
-2 D_{k}^{ii} & i=j 
\end{matrix} \right . 
\qquad D^k_{ij} = \left \{  \begin{matrix}
(-1)^{|e_i||e_j|}C^{k}_{ij} & i \neq j \\ \\
-2 C^{k}_{ii}& i=j 
\end{matrix} \right .\]
\noindent
These will follow directly from the definitions of linear duality:
\[ \left < [x^*,y^*]_{\g^*}, z \right > = \left < x^* \otimes y^*,\delta_{\g}(z) \right > 
\quad 
\left <\delta_{\g^*}(z^*), x \otimes y \right > = \left <z^*, [x,y]_{\g} \right> \]
\noindent
where we assume $x,y,z \in \g$ and $x^*, y^*, z^* \in \g^*$ are homogeneous. Clearly $C^{k}_{ii} =
D_{k}^{ii} = 0$ unless $e_i$ (and hence $e_i^*$) is odd. 

In this setup, the opposite super Lie bialgebra structure on $\g^*$ can be defined as follows
\footnotemark :
\[ [e_i^*,e_j^*]_{(\g^*)^{op}} = [e_j^*,e_i^*]_{\g^*} \qquad  \delta_{(\g^*)^{op}} (e_k^*) =
\delta_{\g^*} (e_k^*) \]
\noindent
or equivalently:
\[ [\cdot,\cdot]_{(\g^*)^{op}} = (-\delta_{\g})^* \qquad  \delta_{(\g^*)^{op}} = \delta_{\g^*}  \]
\noindent
Note that we are only taking the opposite in terms of the Lie superalgebra structure. 

\footnotetext{To compare with \cite{Gou}, note that $T_s(a \wedge b) = -(a \wedge b)$. Here $T_s$ is
defined on homogeneous elements by $T_s(a\otimes b) = (-1)^{|a||b|}b \otimes a$. In
other words, it is the permutation map in the category of super vector spaces.}

The Drinfeld double $\d$ of $\g$ will be defined as the super Lie bialgebra with the underlying
graded vector space identified with $\g \oplus \g^* \cong \g \oplus (\g^*)^{op}$. In order to
define a Lie superalgebra structure on $\d$, we first define a non-degenerate inner product $\left
< \cdot , \cdot \right >$ on $\d$ by asserting super-symmetry: 
\[ \left < x^*, y \right > = (-1)^{|x^*||y|} \left < y, x^* \right > \]
\noindent
and the isotropy of the subspaces $\g$ and $\g^*$:
\[ \left <\g ,\g \right > = \left < \g^*, \g^* \right > = 0 \]
\noindent
(This choice of notation is intentional, and is meant to agree with that for the duality). We will
require invariance of this form, which in terms of the bracket on $\d$ translates to: 
\[ \left < [x^*,y]_{\d}, z \right > = \left < x^*, [y, z]_{\d} \right > \qquad  \left < x^*,
[y^*,z]_{\d} \right > = \left < [x^*, y^*]_{\d}, z \right > \]
\noindent
Then the condition that $[\cdot, \cdot]_{\d}$ restricts to $[\cdot,\cdot]_{\g}$ and 
$[\cdot,\cdot]_{(\g^*)^{op}}$, respectively, on $\g$ and $\g^*$ yields the following description of
$[\cdot, \cdot]_{\d}$ in terms of the structure constants of $\g$ and $\g^*$:
\begin{align*}
[e_i,e_j]_{\d} &= \sum_{k} C^{k}_{ij}e_k \\
[e_i^*,e_j^*]_{\d} &= \sum_{k} C^{ji}_{k} e_k^* \\
[e_i^*,e_j]_{\d} &= \sum_{k} C^{ik}_{j}e_k + \sum_{k} C^{i}_{jk} e_k^*
\end{align*}

The super Lie bialgebra structure on $\d$ is defined to make the natural injections $\g \rightarrow
\d$ and $(\g^*)^{op} \rightarrow \d$ embeddings of super Lie bialgebras, and hence is given by:
\[ \delta_{\d} = \delta_{\g} + \delta_{(\g^*)^{op}} = \delta_{\g} + \delta_{\g^*} \]

With the given structures, it can be shown (see \cite{Gou} for details) that $\d$ is a
quasitriangular super Lie bialgebra, with the $r-$matrix:
\[ \sum_i (-1)^{|e_i|} e_i^* \otimes e_i \]

Although the superization of the main concepts we are using may seem straightforward, it can be
shown that several unexpected situations come up during the process. For instance, see \cite{Lei}
for some interesting examples of Manin triples and a few such unexpected phenomena.

\section{Defining the $r-$matrix $r(f)$}
\label{SectionExamplermatrix}

From now on, let $\g = sl(2,1)$. Define:  
\begin{align*} 
& f(E_{11}+E_{33}) = 0 && f(E_{22}+E_{33}) = E_{22}+E_{33} \\
& f(E_{21}) = 0 && f(E_{12}) = E_{12} \\
& f(E_{23}) = 0 && f(E_{13}) = E_{13} \\
& f(E_{31}) = -E_{13} && f(E_{32}) = E_{23} + E_{32}  
\end{align*}
\noindent
and extend $f$ to a linear map on $\g$. It is easy to check that for any $x$, $y \in \g$, this
function satisfies:
\begin{align*}
(f-1)[f(x),f(y)] = f([(f-1)(x),(f-1)(y)])
\label{YBforf}
\end{align*}
\noindent  
which is equivalent to the associated $2-$tensor being an $r-$matrix, (see Lemma $1$ of \cite{Kar}). 

We write the quadratic Casimir element of $\g$  or equivalently the invariant tensor in $\g \otimes
\g$:
\begin{align*}
\Omega \quad &= \quad (E_{11}+ E_{33}) \otimes (-E_{22}-E_{33}) &+& \quad
 (-E_{22}-E_{33}) \otimes (E_{11}+E_{33}) \\
& + \quad ( E_{12}\otimes E_{21} + E_{21}\otimes E_{12} ) &+& \quad
(-E_{13}\otimes E_{31} + E_{31}\otimes E_{13} ) \\
& + \quad
(-E_{23}\otimes E_{32} + E_{32}\otimes E_{23} ). &&
\end{align*}
\noindent
Defining $r(f)$ to be the $2-$tensor $(f \otimes
1)\Omega$, we get: 
\[ r(f) = r_0 + E_{12}\otimes E_{21} - E_{13} \otimes E_{31} +
E_{32} \otimes E_{23}  - E_{13}\otimes E_{13} +   E_{23} \otimes
E_{23}  \]
\noindent
where $r_0 = (-E_{22} - E_{33}) \otimes (E_{11} + E_{33}).$ It is easy to see that $r(f)$ satisfies:
\begin{align}
\label{unitarity}
r + T_s(r) = \Omega. 
\end{align}
\noindent
Recall here that $T_s$ is the permutation map in the category of super vector spaces.

$r(f)$ 
does not allow a straightforward Belavin-Drinfeld type description. In fact we can prove that the
two subsuperalgebras $\Im(f)$ and $\Im(f-1)$ will never be simultaneously isomorphic to root 
subsuperalgebras. The corresponding subsuperalgebras for functions constructible by Belavin-Drinfeld
type data are always root subsuperalgebras. 

All $r-$matrices on a simple Lie algebra satisfying Eqn.\eqref{unitarity} are constructible by
Belavin-Drinfeld type data. Thus, the existence in the super case of an $r-$matrix satisfying this
equation but not allowing a Belavin-Drinfeld type description provides us with yet another example
when the graded case is more involved than the non-graded case. 

\section{The Super Lie Bialgebra Structure Associated to $r(f)$}
\label{SectionLieBialgebraStructure} 

In the rest of this note, we will concentrate on the Lie bialgebra structure on $\g$ associated to
$r(f) = (f \otimes 1)\Omega$, where $f$ is the linear map introduced in Section
\ref{SectionExamplermatrix} above. After explicitly describing this super Lie bialgebra structure we
will compare it with the standard structure.

\subsection{The cocommutator $\delta_f$}
\label{Cocommutatorforf}

We first describe this structure in terms of a cocommutator $\delta_f$, the coboundary $d
(r(f))$ of the $r-$matrix $r(f)$. We have:
\begin{align*} 
r(f) &= (-E_{22} - E_{33}) \otimes (E_{11} + E_{33})
+ E_{12}\otimes E_{21} - E_{13} \otimes E_{31} 
\\
& \qquad +
E_{32} \otimes E_{23}  - E_{13}\otimes E_{13} +   E_{23} \otimes
E_{23}  
\end{align*}
\noindent
To compute $\delta_f$ we use $\delta_f(g) = d(r(f))(g) = g \cdot r(f) = [g \otimes 1 + 1 \otimes g,
r(f)]$, and this gives us: 
\begin{align*}
\delta_f(E_{11}+E_{33}) &= -E_{23} \wedge E_{23} \\ 
\delta_f(E_{22}+E_{33}) &= E_{13}\wedge E_{13} \\
\delta_f(E_{21}) &= E_{21} \wedge (E_{11}+E_{33}) - E_{23} \wedge (E_{13}+E_{31})\\ 
\delta_f(E_{12}) &= E_{12} \wedge (-E_{22}-E_{33}) - (-E_{13}) \wedge (E_{23}+E_{32}) \\
\delta_f(E_{23}) &= 0 \\ 
\delta_f(E_{13}) &= 0 \\
\delta_f(E_{31}) &= (E_{13}+E_{31}) \wedge (E_{11}+E_{33}) + E_{21} \wedge E_{23}\\ 
\delta_f(E_{32}) &= (E_{23}+E_{32}) \wedge  (-E_{22}-E_{33}) + (-E_{12} )\wedge E_{13}
\end{align*}
\noindent
where we use the notation: $a \wedge b = a \otimes b - (-1)^{|a||b|}b \otimes a$ for any two
homogeneous elements $a$, $b$.

\subsection{Two subalgebras of $\g$}
\label{TwoSubalgebras}

Consider the following subspaces defined by $f$:
\begin{align*}
S_1 =& \Im(f-1) = \left < E_{11} + E_{33}, E_{21}, E_{23}, E_{13}+E_{31} \right > 
\\
S_2 =& \Im(f) = \left < E_{22} + E_{33}, E_{12}, E_{13}, E_{23}+E_{32} \right >
\end{align*}
\noindent
The fact that $r(f)$ is an $r-$matrix satisfying Eqn.\eqref{unitarity} implies that these image
subspaces are indeed Lie subsuperalgebras of $\g$, (see Lemma 4 of \cite{Kar}). In fact it is not
difficult to see that both $S_i$ are isomorphic as Lie superalgebras to a four-dimensional Lie
superalgebra:
\[ \s =  \s_{\overline{0}} \oplus \s_{\overline{1}}; \qquad \qquad \s_{\overline{0}} = \left < h, x
\right >\quad \s_{\overline{1}} = \left < y_1, y_2 \right > \]
\noindent
with the following relations:
\begin{align*}
[h, x] = -x, & &[h, y_1] = -y_1, & &
[x, y_2] = y_1,  \\
[y_1,y_2] = x, & &
[y_2, y_2] = 2h,
\end{align*}
\noindent
(any other commutator will be equal to zero). Therefore, we can write:
\[ \g \cong \s \oplus \s \]
\noindent
where the direct sum is the direct sum of graded vector spaces. Note that $\s$ is solvable. 

Next we compute the restriction of $\delta_f$ onto the $S_i$. This is straightforward; on $S_1$ we
get:
\begin{align*}
\delta_f(E_{11}+E_{33}) &= -E_{23} \wedge E_{23} \\ 
\delta_f(E_{21}) &= E_{21} \wedge (E_{11}+E_{33}) - E_{23} \wedge (E_{13}+E_{31})\\ 
\delta_f(E_{23}) &= 0 \\ 
\delta_f(E_{13} + E_{31}) &= (E_{13}+E_{31}) \wedge (E_{11}+E_{33}) + E_{21} \wedge E_{23}
\end{align*}
\noindent
and the restriction onto $S_2$ is given by:
\begin{align*}
\delta_f(E_{22}+E_{33}) &= E_{13}\wedge E_{13} \\
\delta_f(E_{12}) &= E_{12} \wedge (-E_{22}-E_{33}) - (-E_{13}) \wedge (E_{23}+E_{32}) \\
\delta_f(E_{13}) &= 0 \\
\delta_f(E_{23} + E_{32}) &= (E_{23}+E_{32}) \wedge  (-E_{22}-E_{33}) + (-E_{12} )\wedge E_{13}
\end{align*}
\noindent
In fact we can see that this gives Lie bialgebra structures to the $S_i$. Hence $(S_i, \delta_f
|_{S_i})$ are actually Lie subbialgebras of $\g$. 

Let us now compute the Lie brackets defined on $S_i^*$ by $\delta_f|_{S_i}$. For simplicity, we will
work with the isomorphic super Lie bialgebras on $\s$. Denote the associated cocommutators on $\s$
by $\delta_i$; in other words, define $\delta_1$ and $\delta_2$ so that $(S_1, \delta_f |_{S_1})$
is isomorphic to the super Lie bialgebra $(\s, \delta_1)$ and $(S_2, \delta_f |_{S_2})$ is
isomorphic to the super Lie bialgebra $(\s, \delta_2)$. Clearly we will have:
\begin{align*}
&\delta_1(h) = -y_1 \wedge y_1 && \delta_2(h) = y_1 \wedge y_1\\ 
&\delta_1(x) = x \wedge h - y_1 \wedge y_2 && \delta_2(x) = -(x \wedge h - y_1 \wedge y_2)\\ 
&\delta_1(y_1) = 0 && \delta_2(y_1) = 0\\ 
&\delta_1(y_2) = y_2 \wedge h + x \wedge y_1 && \delta_2(y_2) = -(y_2 \wedge h + x \wedge y_1)
\end{align*}
\noindent
and we have $\delta_1 = -\delta_2$. In particular, we see that $(S_2, \delta_f |_{S_2})$ is
(isomorphic to) the opposite super Lie bialgebra of $(S_1, \delta_f |_{S_1})$.

Recall that the Lie bracket $[\cdot \; , \cdot ]_1$ on the dual $\s^*$ associated to $\delta_1$ can
uniquely be determined by the following: For any two elements $\alpha$, $\beta \in \s^*$ and any
element $s \in \s$, we have:
\[ \left < [\alpha, \beta]_1, s \right > = \left < \alpha \otimes \beta, \delta_1(s) \right > \]
\noindent
(where $\left < \alpha,s\right > = \alpha(s)$ is the pairing of $\s$ with its dual $\s^*$). For
instance we have:
\begin{align*}
\left < [y_1^*, y_1^*]_1, h \right > & = \left < y_1^* \otimes y_1^*, \delta_1(h) \right > =
- \left < y_1^* \otimes y_1^*,y_1 \wedge y_1 \right > = -2 \left < y_1^* \otimes y_1^*, y_1 \otimes
y_1 \right > \\
& = -2 (-1)^{|y_1||y_1^*|}(y_1^*(y_1))^2 = 2
\end{align*}
\noindent
and $\left < [y_1^*, y_1^*]_1, s \right > = 0$ for any other basis vectors of $\s$. Therefore we
get:
\[ [y_1^*, y_1^*]_1 = 2h^*. \]
\noindent
Similarly we have:
\begin{align*}
\left < [y_1^*, y_2^*]_1, x \right > & = \left < y_1^* \otimes y_2^*, \delta_1(x) \right > =
\left < y_1^* \otimes y_2^*,x \wedge h - y_1 \wedge y_2 \right > = - \left < y_1^* \otimes y_2^*,
y_1 \otimes y_2 \right > \\
& = - (-1)^{|y_1||y_2^*|}y_1^*(y_1)y_2^*(y_2) = 1
\end{align*}
\noindent
and $\left < [y_1^*, y_2^*]_1, s \right > = 0$ for any other basis vectors of $\s$. Therefore we
get:
\[ [y_1^*, y_2^*]_1 = x^*. \]
\noindent
Likewise, we compute the other brackets on $\s^*$. The nonzero brackets are: 
\begin{align*}
[h^*, x^*]_1 = -x^*, & &[h^*, y_2^*]_1 = -y_2^*, & &
[x^*, y_1^*]_1 = y_2^*,  \\
[y_1^*,y_2^*]_1 = x^*, & &
[y_1^*, y_1^*]_1 = 2h^*,
\end{align*}
At this point, it is easy to notice that this is actually isomorphic to the Lie superalgebra $\s$
itself (via the map:
\[ h^* \mapsto h, \qquad x^* \mapsto x, \qquad y_1^* \mapsto y_2, \qquad y_2^* \mapsto y_1). \]

Similar computations on $(S_2, \delta_f |_{S_2})$ show that the super Lie bialgebra $(\s, \delta_2)$
is also self-dual. In particular, the Lie bracket $[\cdot \; , \cdot ]_2$ on the dual $\s^*$
associated to $\delta_2$ is given by:
\begin{align*}
[h^*, x^*]_2 = x^*, & &[h^*, y_2^*]_2 = y_2^*, & &
[x^*, y_1^*]_2 = -y_2^*,  \\
[y_1^*,y_2^*]_2 = -x^*, & &
[y_1^*, y_1^*]_2 = -2h^*,
\end{align*}
and this is isomorphic to the Lie superalgebra $\s$ itself (via the map:
\[ h^* \mapsto -h, \qquad x^* \mapsto x, \qquad y_1^* \mapsto y_2, \qquad y_2^* \mapsto -y_1). \]

We have seen earlier that $(S_2, \delta_f |_{S_2})$ is (isomorphic to) the opposite super Lie
bialgebra of $(S_1, \delta_f |_{S_1})$. Therefore:
\[ \g \cong \s \oplus (\s^*)^{op} \]
\noindent
where the direct sum is that of graded vector spaces. Here we may assume $\s$ is equipped with the
super Lie bialgebra structure given by $\delta_1$ or $\delta_2$, as both are self-dual. 

\subsection{The Drinfeld double of $\s$}
\label{DrinfeldDoubleofs}

Let $\d$ be the Drinfeld double of $(\s,\delta_2)$, (we can carry out the following using $(\s,
\delta_1)$ instead, our choice is in fact arbitrary). The computations of \ref{TwoSubalgebras} can
be
used to conclude that $\d \cong \g$ as a Lie superalgebra. Explicitly the map $i_1 : \s \rightarrow
\g$ given on the generators of $\s$ as:
\begin{align*}
&h \mapsto E_{22}+E_{33} &&
x \mapsto E_{12} \\
&y_1 \mapsto E_{13} &&
y_2 \mapsto E_{23}+E_{32}
\end{align*}
\noindent
and the map $i_2 :(\s^*)^{op} \rightarrow \g$ given on the generators as:
\begin{align*}
&h^* \mapsto -(E_{11}+E_{33}) &&
x^* \mapsto E_{21} \\
&y_1^* \mapsto -E_{13} - E_{31} &&
y_2^* \mapsto E_{23}
\end{align*}
\noindent
are both super Lie bialgebra homomorphisms, and $\Im(i_1) = S_2$ and $\Im(i_2) = S_1$.
 
The inner product defined on $\d$ is given by:
\[ \left < s_1 + \alpha_1, s_2+\alpha_2 \right > =_{def} \alpha_1(s_2) +
(-1)^{|\alpha_2||s_1|}\alpha_2(s_1) \]
\noindent
Clearly $S_i$ are both isotropic with respect to this form. We only need to consider
$\left <\alpha_1,s_2 \right >$ and $\left <s_1,\alpha_2 \right >$ where $s_i \in \Im(i_1)$ and
$\alpha_i \in \Im(i_2)$. Now,
\[ \left <\alpha_1,s_2 \right > = \alpha_1(s_2), \qquad \left <s_1,\alpha_2 \right > =
(-1)^{|\alpha_2||s_1|}\alpha_2(s_1) \]
\noindent
and we see that this form actually coincides with the super trace form on $\g$. For example we can
compute:
\begin{align*}
\left <E_{23}, E_{32}\right > &= \left <E_{23}, (E_{23} + E_{32}) - E_{23}\right > = \left <E_{23},
(E_{23} + E_{32})\right > \\
& = \left <i_2(y_2^*), i_1(y_2)\right > = y_2^*(y_2) = 1 \\
\left <E_{13}, E_{31}\right > &= \left <E_{13}, -(-E_{13} - E_{31}) - E_{13}\right > =  \left
<E_{13}, -(-E_{13}- E_{31})\right > \\
&= \left <i_1(y_1), -i_2(y_1^*)\right > = -(-1)^{|y_1||y_1^*|}y_1^*(y_1) = 1 
\end{align*}

We know the super Lie algebra structure when restricted to $S_1$ and $S_2$. From the invariance of
the form we can find the mixed brackets. In other words we use: 
\begin{align*}
\left< [\alpha_1,s_1],s_2 \right > = \left < \alpha_1,[s_1,s_2] \right > & & \left < [s_1,
\alpha_1],\alpha_2 \right > = \left < s_1,
[\alpha_1,\alpha_2] \right >
\end{align*}
where $\alpha_i \in \Im(i_2) = S_1$ and $s_i \in \Im(i_1) = S_2$. Some more computation shows that
indeed the Lie superalgebra structure on $\d$ is the usual one, in other words, we find that $\d
\cong \g$ as a Lie superalgebra. 

The super Lie bialgebra structure on $\d$ is given by $\delta|_{\Im(\s)} + \delta|_{\Im(\s^*)^{op}}$
or equivalently by $\delta|_{\Im(\s)} + \delta|_{\Im(\s^*)}$. But this is equal to $\delta_f|_{S_1}
+ \delta_f|_{S_2} = \delta_f$. Hence, $\g = sl(2,1)$ with the super Lie bialgebra structure
$\delta_f$ is the Drinfeld double of the four dimensional solvable Lie superalgebra $\s$. 

\subsection{The standard structure}
\label{StandardStructure}

It may be interesting to compare the above with the standard Lie bialgebra structure on $\g$. Note
that
\[ r_s = (-E_{22} - E_{33}) \otimes (E_{11} + E_{33})
+ E_{12}\otimes E_{21} - E_{13} \otimes E_{31} +
E_{32} \otimes E_{23}  \]
\noindent
gives us an $r-$matrix satisfying Eqn.\eqref{unitarity} and hence is a standard $r-$matrix
constructed by Belavin-Drinfeld type data. Using $\delta_s = d(r_s)$ to find the associated
cocommutator $\delta_s$ we see that:
\begin{align*}
\delta_s(E_{11}+E_{33}) &= 0 \\ 
\delta_s(E_{22}+E_{33}) &= 0 \\
\delta_s(E_{21}) &= E_{21} \wedge (E_{11}+E_{33}) - E_{23} \wedge E_{31}\\ 
\delta_s(E_{12}) &= E_{12} \wedge (-E_{22}-E_{33}) - (-E_{13}) \wedge E_{32} \\
\delta_s(E_{23}) &= 0 \\ 
\delta_s(E_{13}) &= 0 \\
\delta_s(E_{31}) &= E_{31} \wedge (E_{11}+E_{33}) \\ 
\delta_s(E_{32}) &= E_{32} \wedge  (-E_{22}-E_{33}) 
\end{align*}
\noindent
Restriction of $\delta_s$ onto the $S_i$ would not give us well-defined maps on the $S_i$. Instead,
we may consider its restriction onto two other solvable subalgebras of $\g$: Define 
\begin{align*}
T_1 =& \left < E_{11} + E_{33}, E_{21}, E_{23}, E_{31} \right > ,
\\
T_2 =& \left < E_{22} + E_{33}, E_{12}, E_{13}, E_{32} \right >.
\end{align*}
\noindent
Then on $T_1$ we get:
\begin{align*}
\delta_s(E_{11}+E_{33}) &= 0 \\ 
\delta_s(E_{21}) &= E_{21} \wedge (E_{11}+E_{33}) - E_{23} \wedge E_{31}\\ 
\delta_s(E_{23}) &= 0 \\ 
\delta_s(E_{31}) &= E_{31} \wedge (E_{11}+E_{33}) 
\end{align*}
\noindent
and the restriction of $\delta_s$ onto $T_2$ is given by:
\begin{align*}
\delta_s(E_{22}+E_{33}) &= 0 \\
\delta_s(E_{12}) &= E_{12} \wedge (-E_{22}-E_{33}) - (-E_{13}) \wedge E_{32} \\
\delta_s(E_{13}) &= 0 \\
\delta_s(E_{32}) &= E_{32} \wedge  (-E_{22}-E_{33}) 
\end{align*}
\noindent
It is easy to see that $(T_i, \delta_s|_{T_i})$ is a Lie subbialgebra of $(\g,\delta_s)$ for each
$i$. 

It is clear that $T_i \cong \t$ where $\t$ is the four dimensional solvable Lie superalgebra
\[ \t = \t_{\overline{0}} \oplus \t_{\overline{1}}; \qquad \qquad \t_{\overline{0}} = \left < h, x
\right >\quad \t_{\overline{1}} = \left < y_1, y_2 \right > \]
\noindent
with the following relations:
\begin{align*}
[h, x] = -x, & &[h, y_1] = -y_1, && [y_1,y_2] = x,
\end{align*}
\noindent
(any other commutator will be equal to zero). Therefore, we can write:
\[ \g \cong \t \oplus \t \]
\noindent
where the direct sum is that of graded vector spaces.

Let us now compute the Lie brackets defined on $T_i^*$ by $\delta_s|_{T_i}$. For simplicity, we will
work with the isomorphic super Lie bialgebras on $\t$. Denote the associated cocommutators on $\t$
by ${\delta_s}_i$; in other words, define ${\delta_s}_1$ and ${\delta_s}_2$ so that $(T_1, \delta_s
|_{T_1})$ is isomorphic to the super Lie bialgebra $(\t, {\delta_s}_1)$ and $(T_2, \delta_s
|_{T_2})$ is isomorphic to the super Lie bialgebra $(\t, {\delta_s}_2)$. Clearly we will have:
\begin{align*}
&{\delta_s}_1(h) = 0 && {\delta_s}_2(h) = 0\\ 
&{\delta_s}_1(x) = x \wedge h - y_1 \wedge y_2 && {\delta_s}_2(x) = -(x \wedge h - y_1 \wedge y_2)\\ 
&{\delta_s}_1(y_1) = 0 && {\delta_s}_2(y_1) = 0\\ 
&{\delta_s}_1(y_2) = y_2 \wedge h  && {\delta_s}_2(y_2) = -y_2 \wedge h 
\end{align*}
\noindent
and we have ${\delta_s}_1 = -{\delta_s}_2$. In particular, we see that $(T_2, \delta_s |_{T_2})$ is
(isomorphic to) the opposite super Lie bialgebra of $(T_1, \delta_s |_{T_1})$.

At this point, the Lie bracket $[\cdot \; , \cdot ]_1$ on the dual $\t^*$ associated to
${\delta_s}_1$ can be easily determined. The nonzero brackets are: 
\begin{align*}
[h^*, x^*]_1 = -x^*, & &[h^*, y_2^*]_1 = -y_2^*, & &
[y_1^*,y_2^*]_1 = x^*, 
\end{align*}
Notice that this is actually isomorphic to the Lie superalgebra $\t$ itself (via the map:
\[ h^* \mapsto h, \qquad x^* \mapsto x, \qquad y_1^* \mapsto y_2, \qquad y_2^* \mapsto y_1). \]

As $(T_2, \delta_s |_{T_2})$ is (isomorphic to) the opposite super Lie bialgebra of $(T_1, \delta_s
|_{T_1})$, we get:
\[ \g \cong \t \oplus (\t^*)^{op} \]
\noindent
where the direct sum is that of graded vector spaces, and we are considering $\t$ with the super Lie
bialgebra structure given by ${\delta_s}_1$. Of course, we can show that $(\t, {\delta_s}_2)$ is
self-dual as well, and hence the above identity would still hold if we assumed that the super Lie
bialgebra structure on $\t$ is the one asociated to ${\delta_s}_2$.

Now arguments similar to those in \ref{TwoSubalgebras} and \ref{DrinfeldDoubleofs}  can be used to
conclude that $\g$ is isomorphic to the double of $\t$ as a Lie superalgebra. Explicitly the map
${i_s}_1 : \t \rightarrow \g$ given on the generators of $\t$ as:
\begin{align*}
&h \mapsto E_{22}+E_{33} &&
x \mapsto E_{12} \\
&y_1 \mapsto E_{13} &&
y_2 \mapsto E_{32}
\end{align*}
\noindent
and the map ${i_s}_2 :(\t^*)^{op} \rightarrow \g$ given on the generators as:
\begin{align*}
&h^* \mapsto -(E_{11}+E_{33}) &&
x^* \mapsto E_{21} \\
&y_1^* \mapsto -E_{31} &&
y_2^* \mapsto E_{23}
\end{align*}
\noindent
are both super Lie bialgebra homomorphisms, and $\Im({i_s}_1) = T_2$ and $\Im({i_s}_2) = T_1$.
 
The inner product defined on $\g$ by this double structure is given by:
\[ \left < s_1 + \alpha_1, s_2+\alpha_2 \right > =_{def} \alpha_1(s_2) +
(-1)^{|\alpha_2||s_1|}\alpha_2(s_1) \]
\noindent
Clearly $T_i$ are both isotropic with respect to this form. We only need to consider
$\left < \alpha_1,t_2 \right >$ and $\left < t_1,\alpha_2 \right >$ where $t_i \in \Im({i_s}_1)$ and
$\alpha_i \in \Im({i_s}_2)$,
but
\[ \left < \alpha_1,s_2 \right > = \alpha_1(s_2), \qquad \left < s_1,\alpha_2 \right > =
(-1)^{|\alpha_2||s_1|}\alpha_2(s_1) \]
\noindent
and we see that this form also coincides with the super trace form on $\g$. 

We know the super Lie algebra structure when restricted to $T_1$ and $T_2$. From the invariance of
the form we can find the mixed brackets. In other words we use: 
\begin{align*}
\left <[\alpha_1,t_1],t_2 \right > = \left < \alpha_1,[t_1,t_2] \right > & & \left <[t_1,
\alpha_1],\alpha_2 \right > = \left <t_1,
[\alpha_1,\alpha_2] \right > 
\end{align*}
where $\alpha_i \in \Im({i_s}_2) = T_1$ and $t_i \in \Im({i_s}_1) = T_2$. Some more computation
shows that indeed the Lie superalgebra structure on $\g$ is the usual one. 

The super Lie bialgebra structure on $\g$ coming from this double construction is given by
$\delta|_{\Im(\t)} + \delta|_{\Im(\t^*)^{op}}$ or equivalently by $\delta|_{\Im(\t)} +
\delta|_{\Im(\t^*)}$. But this is equal to $\delta_s|_{T_1} + \delta_s|_{T_2} = \delta_s$.
Therefore, $\g = sl(2,1)$ with the standard super Lie bialgebra structure is the Drinfeld double of
the four dimensional solvable Lie superalgebra $\t$ (equipped with the super Lie bialgebra structure
given either by ${\delta_s}_1$ or ${\delta_s}_2$). 

\section{Conclusion}
\label{SectionConclusion}

Unlike in the non-graded case, super $r-$matrices satisfying Eqn.\eqref{unitarity} may not be
obtained via a simple modification of the Belavin-Drinfeld construction. In this note, we have
studied the Lie bialgebra structure associated to one such $r-$matrix on $sl(2,1)$, and we have
shown that it has a nice description as the double of a four dimensional subalgebra. In the
non-graded case, such structures only arise from twists of the standard $r-$matrix. Our $r-$matrix
is not of this form, but shows similarities to such. These similarities may lead to an
understanding of these special types of $r-$matrices that do not fit a Belavin-Drinfeld type
description.

\end{document}